\documentclass{amsart}

\theoremstyle{plain}

\newtheorem{thm}{Theorem}

\newtheorem{lemma}[thm]{Lemma}

\usepackage{amscd,amsfonts,amssymb,latexsym}

\theoremstyle{remark}

\theoremstyle{definition}
\newtheorem{defn}[thm]{Definition}

\newtheorem{example}[thm]{Example}

\newtheorem*{notation}{Notation}
\newtheorem*{construction}{Construction}

\newcommand{\D}{\Delta}

\renewcommand{\L}{\Lambda}

\newcommand{\Pic}{\operatorname{Pic}}

\newcommand{\Picent}{\operatorname{Picent}}
\newcommand{\Aut}{\operatorname{Aut}}
\newcommand{\Inn}{\operatorname{Inn}}
\newcommand{\rad}{\operatorname{rad}}

\newcommand{\End}{\operatorname{End}}

\newcommand{\m}{\mathfrak{m}}
\newcommand{\p}{\mathfrak{p}}
\newcommand{\q}{\mathfrak{q}}
\newcommand{\e}{\varepsilon}

\newcommand{\Z}{\ensuremath{\mathbb{Z}}}

\newcommand{\Q}{\ensuremath{\mathbb{Q}}}
\begin{document}
\title{Strongly Graded Hereditary Orders}
\author{Jeremy Haefner}
\address{Department of Mathematics, University of Colorado, Colorado Springs, CO 80933-7150}
\email{haefner@math.uccs.edu}
\author{Christopher J. Pappacena}
\address{Department of Mathematics, Baylor University, Waco, TX 76798}
\email{Chris\_\;Pappacena@baylor.edu}
\dedicatory{This paper is dedicated to the memory of Professor Dennis Estes.}
\date{\today}
\thanks{This research was funded in part by a Baylor University Research Grant.}

\begin{abstract} Let $R$ be a Dedekind domain with global quotient field $K$.
The purpose of this note is to provide a characterization of when a strongly
graded $R$-order with semiprime $1$-component is hereditary. This generalizes
earlier work by the first author and G. Janusz in (J. Haefner and G. Janusz,
{\it Hereditary crossed products}, Trans. Amer. Math. Soc. \textbf{352} (2000),
3381-3410). \end{abstract}

\maketitle

Recall that, for a Dedekind domain $R$ with quotient field $K$, an
$R$-\emph{order} in a separable $K$-algebra $A$ is a module-finite $R$-algebra
$\Lambda$, contained in $A$, such that  $K\Lambda=A$. For a group $G$, we say
that the $R$-order $\Lambda$ is {\it strongly $G$-graded} provided there is a
decomposition $\Lambda = \oplus_{g\in G} \Lambda_g$ with
$\Lambda_g\Lambda_h =\Lambda_{gh}$ for all $g, h \in G$. If $1$ denotes the
identity element of $G$, then $\Lambda_1$ is a subring of $\Lambda$, which we
denote by $\Delta$. We write $\Lambda = \Delta(G)$ to indicate
that $\Lambda$ is strongly $G$-graded with identity component $\Delta$. In this
note, we consider the following problem:

\begin{quote}
\noindent{\bf The hereditary problem for
strongly graded orders:} Determine necessary and sufficient conditions on $G$,
the grading imposed by $G$, and $\Delta$ to ensure that $\Lambda$ is hereditary.
\end{quote}

Our general solution to this problem appears in Theorem
\ref{main thm}. The idea of the proof is to use Morita theory to reduce to the
case where $\L$ is a crossed product and then apply a result of \cite{HJ}, which
we describe next.

Recall that a strongly $G$-graded $R$-order $\L=\Delta(G)$ is a {\it crossed
product order} provided for each $g\in G$, $\L_g \cong \Delta$ as left
$\Delta$-modules.  In this case, there exist $u_g\in \Lambda^*$ (the unit group
of $\L$) such that $\L_g=\Delta u_g$ for all $g\in G$.  Moreover, there exist a group
homomorphism $\alpha: G \to \Aut_R(\Delta)$ (the ``action of $G$ on $\Delta$")
and a cocycle $\tau \in Z^2(G, R^*)$ (the ``twisting of the action of $G$") such
that the multiplication in $\L$ is given by $u_g\delta=\alpha(g)(\delta)u_g$ for
$\delta\in\Delta$,  and $u_gu_h=\tau(g,h)u_{gh}$.   (See \cite{P} for more
details on this construction.)

If $\Delta(G)$ is a crossed product order with action
$\alpha$, we say that a subgroup $H$ of $G$ {\it acts as central outer
automorphims of $\Delta$} provided $\alpha(H)\cap\Inn(\Delta)=1$. The main
result of \cite{HJ} is that, if  $\L=\Delta(G)$ is a crossed product order, then $\L$
is hereditary if and only if $\Delta$ is hereditary and, for each maximal ideal
$\m$ of $R$ containing a prime divisor $p$ of $|G|$, any $p$-Sylow subgroup of
$G$ acts as central outer automorphisms of $\hat\Delta_\m$.

\begin{defn}
Given a strongly $G$-graded ring $\L=\Delta(G)$ and $g\in G$, we say $g$ is {\it
inner on $\Delta$} provided $\L_g \cong \Delta$ as $\Delta$-bimodules.
Otherwise, $g$ is {\it outer on $\Delta$}. For a subgroup $H$ of $G$, set
$$\Inn_{\Delta}(H)=\{h\in H : h \hbox{ is inner on $\Delta$} \}$$ We say that
$H$ is {\it inner on $\Delta$} or {\it inner grades $\Delta$} if $\Inn_{\Delta}(H)=H$ and it is {\it outer on
$\Delta$} or {\it outer grades $\Delta$} provided $\Inn_{\Delta}(H) = 1$. \end{defn}

We remark that the above
definitions do indeed generalize the classical notion of inner actions. To see
this, suppose $\L = \Delta(G)$ is a crossed product order with group action $\alpha: G \to \Aut_R(\Delta)$ such that $\alpha(g) \in \Inn_R(\Delta)$ for some $g \in G$. Then it is easy to see  that  $\Lambda_g=\Delta u_g \cong \Delta$ as bimodules. Hence, if $g$
acts as an inner automorphism on $\Delta$ in the classical sense, it is inner on
$\Delta$ in the sense defined above.

\begin{lemma}  Let
$\L=\Delta(G)$ be strongly graded by $G$. Then, for any subgroup $H$ of $G$,
$\normalfont\Inn_\Delta(g^{-1}Hg)=g^{-1}\Inn_\Delta(H)g$. \label{conj lemma} \end{lemma}

\begin{proof} Suppose that $h\in\Inn_\Delta(H)$.  Then $\Delta_h\cong \Delta_1$
as bimodules.  It follows that $\Delta_{g^{-1}}\Delta_h\Delta_g\cong \Delta_1$
as bimodules as well.  This shows that $g^{-1}\Inn_\Delta(H)g\subseteq
\Inn_\Delta(g^{-1}Hg)$.  For the converse, let $\Delta_{g^{-1}hg}\cong \Delta_1$
as bimodules.  Since $\Delta$ is strongly graded,
$\Delta_{g^{-1}hg}=\Delta_{g^{-1}}\Delta_h\Delta_g$.  It follows that
$\Delta_h\cong \Delta_1$ as bimodules, proving the reverse inclusion.
\end{proof}

\begin{lemma}  Suppose that $R$ is a complete DVR, and that
$\L=\Delta(G)$ is a strongly-graded $R$-order such that $\Delta$ is prime
and basic.  Then $\Delta(G)$ is a crossed product order. \label{basic} \end{lemma}

\begin{proof} Since $R$ is a complete DVR, $\Delta$ is semiperfect.  Consequently there exist only finitely many indecomposable projective left $\Delta$-modules up to isomorphism. Let $\{P_1,\dots,P_t\}$ be representatives for these isomorphism classes. Then there exist positive integers $m_i$ such that $\Delta\cong\oplus_{i=1}^t P_i^{(m_i)}$ as left $\Delta$-modules.

 For $g\in G$, we have that $\L_g\in \Pic(\Delta)$, because $\L$ is strongly graded.  This implies that $\L_g$, viewed as a left $\D$-module, is isomorphic to a direct sum of the $P_i$.  So, we may write $\L_g\cong\oplus_{i=1}^tP_i^{(n_i)}$, where the $n_i$ are \emph{a priori} nonnegative integers.   Now, since $\L_g$ is a progenerator and $\End(\L_g)\cong \Delta$ as rings, a combinatorial argument shows that in fact $m_i=n_i$ for all $i$, so that $\L_g\cong \Delta$ as left $\Delta$-modules.  Thus $\L$ is a crossed product.
\end{proof}

We can now
state and prove the prime case of our solution to the hereditary problem. We use
the following notation  for the remainder of the paper.
$R$ denotes a Dedekind domain whose quotient field $K$ is a global field and $A$
denotes a separable $K$-algebra. For a maximal ideal $\m$ of $R$, $\hat R_\m$
 denotes the completion of $R$ at $\m$. Similarly, $\Delta$  denotes an
$R$-order in $A$, and $\hat \Delta_\m$  denotes the $\hat R_\m$-order $\hat
R_\m\otimes_R \Delta$.

\begin{thm} \label{prime thm}
Assume that $A$ is simple \emph{(}so  $\Delta$ is prime\emph{)}. Then $\L$
is hereditary if and only if $\Delta$ is hereditary and, for each maximal ideal
$\m$ of $R$ containing a prime divisor $p$ of $|G|$, some \emph{(}hence every\emph{)}
$p$-Sylow subgroup of $G$ is outer on $\hat \Delta_\m$. \end{thm}

\begin{proof} First note that the
induction functor $\L \otimes_\Delta -: \Delta-\hbox{mod} \to \L-\hbox{mod}$
is separable \cite[Proposition 2.2]{HJ}.  It follows by \cite[Proposition
2.3]{HJ} that if $\Lambda$ is hereditary, then so is $\Delta$. Thus, for the
remainder of the proof, we assume that $\Delta$ is hereditary. Next observe that
any two $p$-Sylow subgroups of $G$ are conjugate.  Thus, by Lemma \ref{conj
lemma}, it suffices to verify that $\Inn_{\hat\Delta_\m}(P)=1$ for a single
$p$-Sylow subgroup $P$. By \cite[Theorem 40.5]{MO}, $\L$ is hereditary if and
only if $\hat\L_\m$ is hereditary for each maximal ideal $\m$.  We  show that
$\L$ is hereditary if and only if $\hat\L_\m$ is hereditary for those maximal
ideals $\m$ containing a prime divisor of $|G|$.

To see this, assume
$\hat\L_\m$ is hereditary for those maximal ideals $\m$ containing a prime
divisor of $|G|$. Fix an arbitrary maximal ideal $\m$ of $R$.  If it contains no
prime divisors of $|G|$, then $|G|$ is a unit in $\hat R_\m$.  By Proposition
2.2 of \cite{HJ}, we have that $\hat\L_\m$ is a separable extension of
$\hat\Delta_\m$.  Since $\Delta$ is hereditary, so is $\hat\Delta_\m$ and it
follows from   \cite[Proposition 2.3]{HJ}  that $\hat\L_\m$ is hereditary. Thus,
to verify whether $\L$ is hereditary, it suffices to check $\hat\L_\m$ at those
maximal ideals $\m$ containing prime divisors of $|G|$.

Assume that $\m$
contains a prime divisor $p$ of $|G|$, and fix a $p$-Sylow subgroup $P$ of $G$.
For ease of notation, we write $\hat R$ for $\hat R_\m$, etc.  Since
$\hat\Delta$ is prime hereditary, there is an idempotent $e$ of $\hat\Delta$
such that $e\hat\Delta e$ is basic hereditary. We claim that $e\hat\L e$ is a strongly $G$-graded order with components $e\hat\L_x e$ for $x\in G$. To see this, fix homogeneous
components $e\hat\L_g e, e\hat\L_h e$.  Then,
\begin{eqnarray*}
e\hat\L_g e\cdot e\hat\L_h e&=&e\hat\L_ge\hat\L_he\\
&=&e\hat\L_g\hat\L_1e\hat\L_1\hat\L_he \ \ \mbox{(since $\hat \L$ is strongly
$G$-graded)}\\ &=&e\hat\L_g\hat\Delta\hat\L_he \ \  \mbox{(since
$\hat\L_1=\hat\Delta$ and $\hat\Delta e\hat\Delta=\hat\Delta$)}\\ &=&e\hat\L_{gh}
e \ \ \mbox{(since $\hat \L$ is strongly $G$-graded)}\end{eqnarray*}
which shows that $e\hat\L e$ is strongly graded.

 The orders  $\hat\L$
and $e\hat\L e$ are Morita equivalent via the pair of graded progenerators $\hat\L
e $ and $e\hat\L$. In addition, this equivalence preserves, in a certain sense, the grading of the two orders. This is an example of what is called a {\it graded equivalence}; see \cite{HG} for more information on graded
equivalences.

By Morita
equivalence,  $\hat\L$ is hereditary if and only if $e\hat\L e$ is hereditary.
Now, the identity component of $e\hat\L e$ is $e\hat\Delta e$, which is basic,
prime and hereditary.  Since $e\hat\L e$ is strongly $G$-graded, we see by Lemma \ref{basic} that $e\hat\L e$ is a crossed product order.
Hence, we may apply \cite[Theorem 6.8]{HJ} to conclude that
$e\hat \L e$ is hereditary if and only if $P$ acts as central outer
automorphisms.  As we have remarked above, this is equivalent to saying that $P$
outer grades $e\hat\L e$. To finish the proof, we note that $P$ outer grades $e\hat\L e$ if and only if $P$
outer grades $\hat \L$, because $\hat\L_x \cong \hat\Delta$ as bimodules if and only if  $e\hat\L_x e \cong e\hat\Delta e$ for any $x \in P$. (This uses the fact that $\hat\L e$ and $e\hat\L$ induce a graded equivalence between $\hat\L$ and $e\hat\L e$.) Thus, $\hat\L$ is
hereditary if and only if $\Inn_{\hat\L}(P)=1$. \end{proof}

The above proof
requires us to reduce to the case when $\Delta$ is basic so that $\L$ is a
crossed product order. We present an example that shows that, even over a complete
DVR, a strongly graded order with non-basic hereditary $1$-component need not be
a crossed product order. Thus Lemma \ref{basic} is the best possible, and the proof of
Theorem \ref{prime thm} cannot be simplified in this regard. The example depends
on the following basic construction technique for strongly graded rings, which
is a special case of \cite[p. 23]{NV}.

\begin{construction} Let $\Delta$ be a
ring, and let $X\in\Pic(\Delta)$ have finite order $n$.  Let $G$ denote the
cyclic group of order $n$ with generator $g$.  We construct a
$G$-strongly graded ring $\L$ with $1$-component $\Delta$ as follows: Set
$\L=\oplus_{i=1}^n \L_{g^i}$, where $\L_{g^i}=X^i$, and define the
multiplication by the tensor product. In other words, fix isomorphisms
$X^i\otimes X^j\cong X^{i+j}$ which are compatible in the obvious sense.  Then,
given homogeneous elements $x_i, x_j$ in $\L_{g^i}$, $\L_{g^j}$, respectively,
we define $x_i\cdot x_j=x_i\otimes x_j\in X^{i+j}$ (using the fixed
isomorphism). Note that $\L_{1}=\L_{g^0}=X^0\cong \Delta$, so the $1$-component
is $\Delta$, as claimed. Note also the grading is strong, as $X^i\otimes
X^j\cong X^{i+j}$ for all $i,j$ by construction. 
\end{construction}

\begin{example}[Strongly
graded orders need not be crossed products]   Let $\Delta$ be a non-basic
hereditary order over a complete DVR, and let $\Gamma$ be the associated basic
order $e\Delta e$.  Now, $\Picent(\Delta)\cong\Picent(\Gamma)$ is cyclic (say of order $n$)
generated by $\rad(\Delta)$ (respectively $\rad(\Gamma)$).  The fact that it is
cyclic is in \cite{HJ}, and the fact that it is generated by the radical follows
from \cite{MO}, where it is shown that $\rad(\Delta)$ has the correct order.
Now, if $\Delta$ is not basic, then $\rad(\Delta)^k$ is not principal as a left
ideal for any $1\leq k<n$, and so $\rad(\Delta)^k\not\cong \Delta$ as bimodules for any
$1\leq k<n$.  Thus, the $\Z/n\Z$-strongly graded order $\Delta(G)=\oplus_{k=0}^{n-1} \rad(\Delta)^k$  is
not a crossed product order, even though $e\Delta(G)e$ is.\qed
\label{nonbasic} \end{example}

Having dealt with the prime case, we turn our attention to the
semiprime case. To begin, we investigate the action of $G$ that is imposed on
the central idempotents of $\Delta$.  We fix the following notation.

\begin{notation}  Given
$\L=\Delta(G)$, suppose that
$\Delta=\Delta_1\oplus\cdots\oplus\Delta_t$ is a direct sum of prime rings.  Let
$e_1,\dots,e_t$ denote the orthogonal central idempotents of $\Delta$.  Then the
group $G$ acts on the $e_i$, by $e_i\L_g=\L_g e_{g(i)}$.  (This is a special
case of the action of $\Pic(\Delta)$ on the center of $\Delta$; see \cite[\S
55]{CR2}.) Suppose that this action partitions $\{e_i\}$ into $m$ orbits, and
let $\e_1,\dots,\e_m$ be a set of representatives for the equivalence classes
under this action.  Finally, let $G_i$ denote the stabilizer of $\e_i$.
\end{notation}

\begin{lemma} Assume the above notation.
\begin{enumerate}
\item[(a)] $\L$ is Morita equivalent to $\oplus_{i=1}^m \e_i\L \e_i$.
\item[(b)]
Each $\e_i\L\e_i$ is strongly graded by $G_i$, with identity component
$\e_i\Delta\e_i$, a prime ring.
\end{enumerate} \label{decomp} \end{lemma}

\begin{proof} This is proven in \cite[Theorem 5.4]{H},  under the assumption
that $\Delta$ is maximal.  However, examining the proof, we see that the above
is true without this hypothesis.  \end{proof}

\begin{thm} Let $R$ be a Dedekind
domain with global quotient field $K$. Let $A$ be a semisimple $K$-algebra, and
$\Delta$ be an $R$-order in $A$. \emph{(}Note that $\Delta$ is necessarily
semiprime.\emph{)} Suppose that $\L=\Delta(G)$ is a strongly graded $R$-order.
 Then $\Delta(G)$ is hereditary if and only if $\Delta$ is hereditary, and, in the notation of Lemma \ref{decomp}, the
following conditions hold, for $1\leq i\leq m$.

\begin{quotation}
\noindent For each maximal ideal $\m$ of $R$
containing a prime divisor $p$ of $|G_i|$,  $\normalfont\Inn_{\e_i\hat\Delta\e_i}(P)=1$ for
some \emph{(}hence every\emph{)} $p$-Sylow subgroup $P$ of $G_i$.
\end{quotation}
\label{main thm}
\end{thm}

\begin{proof}  Note that, as in the proof of Theorem \ref{prime
thm}, we conclude that $\Delta$ being hereditary is a necessary condition.
Thus, we may assume $\Delta$ is hereditary, and then $\Delta$ decomposes as a
direct sum of prime rings.  Hence, Lemma \ref{decomp} applies.

Since $\L$
is Morita equivalent to $\oplus_{i=1}^m\e_i\L\e_i$, it follows that $\L$ is hereditary
if and only if each $\e_i\L\e_i$ is.  Now, each $\e_i\L\e_i$ has prime identity
component, so that Theorem \ref{prime thm} applies. \end{proof}

We close this paper with some remarks and examples.
First, the decomposition $\oplus_i\e_i\L\e_i$ depends upon the choice of the
representatives $\e_i$ of the orbits of $\{e_1,\dots,e_t\}$.  If we chose
different representatives $\e_i'$, then the new stabilizer groups $G_i'$ would
be conjugate to the original groups $G_i$.  If we choose $p$-Sylow subgroups
$P$, $P'$ of $G_i$, $G_i'$ respectively, then $\Inn_{\e_i\hat\Delta\e_i}(P)=1$
if and only if $\Inn_{\e_i'\hat\Delta\e_i'}(P')=1$, by Lemma \ref{conj lemma}.
Thus, the choice of representatives does not affect the application of Theorem
\ref{main thm}.

Second, the statement of Theorem \ref{main thm} requires
checking whether or not $\Inn_{\hat\Delta}(P)=1$ at various completions.  It is
not enough to assume that $\Inn_\Delta(P)=1$, i.e. that $P$ outer grades
$\Delta$ globally, because the property of being outer is not a local-global
property.  The next example illustrates this fact.

\begin{example}[Outer grading
is not a local-global property] Let $R$ be the ring of Gaussian integers
$\Z[i]$, and let $K$ denote the quotient field $\Q(i)$.  The prime integer $5$
is contained in exactly two ideals of $R$: $\p=(1+2i)$ and $\q=(1-2i)$.  Let $I$
denote the ideal $(5)$, and let $\Delta$ denote the tiled order
\[\Delta=\left(\begin{array}{ccccc}
R&R&R&R&R\\I&R&R&R&R\\I&I&R&R&R\\I&I&I&R&R\\I&I&I&I&R\end{array}\right).\] We first
compute $\Picent(\Delta)$.  By Fr\"ohlich's Theorem \cite[Theorem 55.25]{CR2}
, there is an isomorphism \[\Picent(\Delta)\cong
\bigoplus_{\mbox{$\m$ maximal}} \Picent(\hat \Delta_\m).\] (Here we
are using that $R$ is a PID.)  Note that, if $\m\neq \p,\q$, then $\hat
I_\m\cong \hat R_\m$, and so $\hat\Delta_\m=M_5(\hat R_\m)$.  In particular,
$\Picent(\hat\Delta_\m)=0$ if $\m\neq \p,\q$.

If $\m=\p$, then $\hat \Delta_\p$ is the unique
basic hereditary order in $M_5(\hat K)$, and so $\Picent(\hat\Delta_\p)\cong
\Z/5\Z$, by \cite[Theorem 39.18]{MO}. Similarly, $\Picent(\hat\Delta_\q)\cong
\Z/5\Z$.  Thus, $\Picent(\Delta)\cong \Z/5\Z\oplus\Z/5\Z$.

Now, let $X$ be the
bimodule generating  the subgroup $0\oplus\Z/5\Z$ of $\Picent(\Delta)$ (i.e. the
component corresponding to $\hat\Delta_\q$), and form the strongly
$\Z/5\Z$-graded order $\L$, where $\L_{g^i}=X^i$.  Note that, globally, $\L$ is
outer graded (because $X^i\cong \Delta$ as bimodules if and only if  $i=0$).
However, if we pass to the completion at $\p$, then $\hat X_\p\cong \hat
\Delta_\p$ as bimodules, by construction.  Thus, it is possible for a global
outer grading to become inner at the completion. \qed \label{outer}
\end{example}

Finally, Theorem \ref{main thm}  requires checking the condition
on the grading at each $\e_i\L\e_i$ (in the notation of Lemma
\ref{decomp}), rather than simply checking the grading on $\L$.   That is, it is
not enough to verify that $\Inn_{\hat\L}(P)=1$ for $p$-Sylow subgroups of each
$G_i$ (or of $G$).  This is because the Morita equivalence between $\L$ and
$\oplus_{i=1}^m\e_i\L\e_i$ is not a graded equivalence. Hence, the property of being an
outer grading is not preserved under this correspondence. Our last example
illustrates this.

\begin{example}[Passing to $\e_i\L\e_i$ can change the
grading]  Let $\Delta$ be a prime, hereditary order, and let $\Delta^{(d)}$
denote the direct sum of $d$ copies of $\Delta$.  The symmetric group $S_d$ acts
on $\Delta^{(d)}$ by permuting coordinates, so that we may form the crossed
product (with trivial twisting) $\Delta^{(d)}(S_d)$ relative to this action. In
the notation of Lemma \ref{decomp}, $G=S_d$ acts transitively on
$\{e_1,\dots,e_d\}$. Let us fix $\e_1=e_1$, so that $G_1=S_{d-1}$, embedded in
$S_d$ as those permutations that fix the first coordinate.

It is straightforward
to compute that $e_1\Delta^{(d)}(S_d)e_1\cong \Delta S_{d-1}$, the ordinary
group ring of $S_{d-1}$ over $\Delta$.  Note that, for \emph{any} $p$-Sylow
subgroup $P$ of $S_{d-1}$, $\Inn_{\Delta S_{d-1}}(P)=P$. However, if we view $P$
as a subgroup of $S_d$, then $\Inn_{\Delta^{(d)}(S_d)}(P)=1$.  This is because,
for $\pi\in P$, $\Delta^{(d)}(S_d)_\pi\not\cong \Delta^{(d)}$ as bimodules.
(The action on the right is twisted by $\pi$, and $\pi$ is not an inner
automorphism of $\Delta$.)  Thus, it is
necessary when applying Theorem \ref{main thm} to consider the grading on each
$\e_i\L\e_i$, and not simply the grading on $\L$.\qed\label{semiprime}
\end{example}

\end{document}